\documentclass[12pt,reqno]{amsart}
\usepackage{amsmath,amsfonts,amssymb,amsthm,amsxtra,array,amscd,calc}

\usepackage{enumerate}
\usepackage[all]{xy}
\usepackage[english]{babel}

\DeclareMathOperator{\ad}{ad}

\DeclareMathOperator{\Ker}{Ker}

\DeclareMathOperator{\GL}{GL}

\DeclareMathOperator{\Mat}{Mat}

\newcommand{\lra}{\longrightarrow}

\newcommand{\ka}{{\mathcal A}}

\newcommand{\ko}{{\mathcal O}}
\newcommand{\ke}{{\mathcal E}}
\newcommand{\kf}{{\mathcal F}}
\newcommand{\kg}{{\mathcal G}}

\newcommand{\ki}{{\mathcal I}}

\newcommand{\ku}{{\mathcal U}}

\newcommand{\C}{{\mathbb C}}

\renewcommand{\P}{{\mathbb P}}

\begin{document}

\newtheorem{sub}{}[section]
\newtheorem{subsub}{}[sub]
\newtheorem{subsubsub}{}[sub]

\title[Babylonian towers of principal bundles]{A Babylonian tower theorem for 
principal bundles over projective spaces}

\author[Biswas]{I. Biswas}
\address{School of Mathematics, Tata Institute of Fundamental
Research, Homi Bhabha Road, Bombay 400005, India}
\email{indranil@math.tifr.res.in}

\author[Coand\u{a}]{I. Coand\u{a}$^{1}$}
\address{Institute of Mathematics of the Romanian Academy, P.O.~Box
1-764, RO-014700 Bucharest, Romania}
\email{Iustin.Coanda@imar.ro}

\author[Trautmann]{G. Trautmann$^{2}$}
\address{FB Mathematik, Universit\"at Kaiserslautern, 
Postfach 3049, D-67653 Kaiserslautern, Germany}
\email{trm@mathematik.uni-kl.de}

\footnotetext[1]{Partially supported by grant 2-CEx06-11-10/25.07.06 
of the Romanian Ministry of Education and Research and by DFG.}
\footnotetext[2]{Partially supported by DFG Schwerpunktprogramm 1094}

\begin{abstract}
We generalise the variant of the Babylonian tower theorem for vector 
bundles on projective spaces proved by I. Coand\u{a} and G. Trautmann 
(2006) to the case of principal $G$-bundles over projective spaces, 
where $G$ is a linear algebraic group defined over an algebraically closed 
field. In course of the proofs some new insight into the structure 
of such principal $G$-bundles is obtained.
\end{abstract} 

\maketitle
 
Let $G$ be a linear algebraic group defined over 
an algebraically closed field $k$.
A principal $G$-bundle over a projective space
${\P}_n$ is called \textit{split} if it admits a
reduction of structure group to a maximal torus of $G$. Since
a finite dimensional
$T$-module, where $T$ is a torus defined over $k$, splits
into a direct sum of one-dimensional $T$-modules, the adjoint
bundle of a split $G$-bundle decomposes into a direct sum of line 
bundles. When $G$ is reductive, also the converse holds:

\vskip3mm
{\bf Proposition 1:}\quad
\textit{Let $G$ be a reductive linear algebraic group.
Let $\ke$ be a principal $G$-bundle over ${\P}_n$
and $\ad(\ke)$ its adjoint 
bundle. If $\ad(\ke)$ splits as 
a direct sum of line bundles, then $\ke$ is split.}
\vskip3mm

For $k=\C$ this is proved in \cite{BGH}, Theorem 4.3, using arguments 
extracted from Grothendieck's paper \cite{g}.  
A proof in any characteristics is presented in Section \ref{prop1}. 
Using this result and the method from \cite{ct} we will prove the 
following:
\vskip3mm

{\bf Theorem 1:}\quad
\textit{Let $G$ be a linear algebraic group, and let $\ke$ be a 
principal $G$-bundle over ${\P}_n$ with adjoint bundle $\ad(\ke).$ 
Assume that $\ke$ can be extended to a principal 
$G$-bundle over ${\P}_{n+m}$
for some $m > {\Sigma}_{i>0}\dim {\fam0 H}^1(\ad(\ke)(-i))$.\\
If ${\fam0 char}(k)=0$ or if ${\fam0 char}(k)=p>0$ and $G$ is reductive, 
then $\ke$ is split as a principal bundle.}
\vskip3mm


Theorem 1 also holds for arbitrary algebraic groups. This follows from 
the proof of Proposition \ref{zartop} below.
When $k={\mathbb C}$ and 
$G$ is a (finite dimensional) complex Lie group, one can use  
arguments analogous to those below to prove that the adjoint bundle of an 
analytic principal $G$-bundle on $\P_n(\C)$ splits as a 
direct sum of line bundles, if it satisfies the extension assumption in 
Theorem 1.

As a byproduct of the proof of Theorem 1 one gets the following theorem.

\vskip3mm
{\bf Theorem 2:}\quad
\textit{Let $\ke$ be a 
principal $G$-bundle over $\P_n$. 
If ${\fam0 H}^1(\ad(\mathcal E)(-i))=0$ for all
$i>0$, then $\ad(\ke)$ splits as a direct sum of 
line bundles. If, moreover, $G$ is reductive then $\ke$ itself 
is split.}
\vskip3mm

Theorem 2 was proved by  Mohan Kumar \cite{mk} under the
assumption that $G=\text{GL}_r(k)$, and it was proved in 
\cite{b} under the assumption that $k = {\mathbb C}$ with
$G$ reductive. Again, the first assertion of Theorem 2 remains valid when 
$k={\mathbb C}$, $G$ is a (finite dimensional) complex Lie group, and 
$\ke $ is a complex analytic principal $G$-bundle.
\vskip3mm


\section{Some non-abelian cohomology}

The following Proposition enables us to work with Zariski open subsets of
$\P_n$ instead of \'etale covers. As before, $k$ will denote
an algebraically closed field.

\begin{sub}\label{zartop}{\bf Proposition:}\quad
a) Let $G$ be an algebraic group over $k$. Then any principal
$G$-bundle over $\P_n$ is Zariski locally trivial.\\
b) For an abelian variety $A$ over $k$, any algebraic principal $A$-bundle 
over $\P_n$ is trivial.
\end{sub}

Proposition \ref{zartop} will be  proved  in Section \ref{pzartop}. 
One should note, however, that b) is not valid for complex analytic principal
bundles with an abelian variety as the structure group.
\vskip10mm

We use the paper of Frenkel \cite{fr} as a 
reference for basic non-abelian cohomology. Let $X$ be a topological space 
and $\mathcal G$ a sheaf of 
(not necessarily abelian) groups. For $U\subset X$ open, let $e_U$ denote 
the unit element of $\kg(U)$.

One defines, using \v{C}ech 
1-cocycles and their equivalence relation, the {\it first cohomology set}
$\text{H}^1(X,\kg)$. It has a {\it marked element} corresponding to 
the 1-cocycle $(e_X)$ on the open cover $\{ X\} $ of $X$. 

If $c\in \text{H}^1(X,\kg)$ is represented by 
$(g_{ij})\in \text{Z}^1(\ku ,\kg)$ 
for some open cover $\ku$ of 
$X$ then one gets a {\it principal} $\kg$-{\it bundle} 
by gluing the sheaves ${\kg}|U_i$ with 
$g_{ij}\cdot - : \kg|U_{ji}\overset{\sim}{\rightarrow} \kg\vert U_{ij}$. 
In this way, $\text{H}^1(X,\kg)$ 
parametrises the isomorphism classes of principal $\kg$-bundles (locally
trivial with respect to the topology of $X$).

A class $c\in \text{H}^1(X,\kg)$ can also be used to define {\it twists}
of sheaves of groups which are acted on by $\kg$. For that let $\ka$ be any 
other sheaf of groups and assume that there is an action 
$\kg\times \ka\to \ka.$ Then a sheaf $\ka^c$ of groups is defined by gluing
with the isomorphisms
$g_{ij}\cdot - : \ka|U_{ji}\overset{\sim}{\rightarrow} \ka\vert U_{ij}$
of the action.

This twisting is obviously an exact functor on the category of $\kg$-sheaves
of groups.

In particular, a new sheaf of groups ${\kg}^c$ is obtained by gluing 
the sheaves 
$\kg\vert U_i$ with $g_{ij}\cdot -\cdot g_{ij}^{-1} : 
{\kg}\vert U_{ji}\overset{\sim}{\rightarrow} {\kg}\vert U_{ij}$ 
with respect to the action of inner automorphisms.
Let 
${\phi}_i : \kg^c\vert U_i\overset{\sim}{\rightarrow}\kg\vert U_i$ 
be the resulting isomorphisms. 

There exists a 
{\it bijection} 
$\text{H}^1(X,\kg^c)\overset{\sim}{\rightarrow} \text{H}^1(X,\kg)$ 
which is constructed by sending (the class of)
$(f_{ij})\in \text{Z}^1(\ku,\kg^c)$ 
to (the class of)
$({\phi}_i(f_{ij})\cdot g_{ij})\in \text{Z}^1(\ku,\mathcal G)$. 
This bijection sends the marked element of 
$\text{H}^1(X,\kg^c)$ to $c$. 

\par Let, now, $1\rightarrow {\mathcal G}^{\prime}\overset{u}{\rightarrow}
{\mathcal G}\overset{p}{\rightarrow} {\mathcal G}^{\prime \prime}
\rightarrow 1$ be a short exact sequence of sheaves of groups on $X$. This 
means that $p$ is an epimorphism of sheaves and that, for every open subset 
$U\subseteq X$, $u(U)$
maps ${\mathcal G}^{\prime}(U)$ isomorphically onto 
$\Ker p(U)$. In particular, every inner automorphism of 
${\mathcal G}(U)$ induces, via $u(U)$, an automorphism of 
${\mathcal G}^{\prime}(U)$. It follows that, if one twists 
$\kg^{\prime}$ 
as above, one obtains a new sheaf of groups $\kg^{\prime c}$ with an 
exact sequence
$1\rightarrow \kg^{\prime c}\rightarrow \kg^c\rightarrow
\kg^{\prime \prime c}\rightarrow 1$.
\vskip5mm

\begin{sub}\label{idim}{\bf Lemma:}\quad
Under the above hypothesis, there exists a canonical map 
$$
{\fam0 H}^1(X,{\mathcal G}^{\prime c})\rightarrow
{\fam0 H}^1(X,\mathcal G)
$$
sending the marked element of 
${\fam0 H}^1(X,{\mathcal G}^{\prime c})$ to $c$ and whose image is 
${\fam0 H}^1(p)^{-1}({\fam0 H}^1(p)(c))$. 

\begin{proof}
One uses the $\text{H}^1$ part of the cohomology exact sequence associated 
to the short exact sequence of sheaves of groups :
$1\rightarrow \kg^{\prime c}\rightarrow \kg^c\rightarrow
\kg^{\prime \prime c}\rightarrow 1$.
\end{proof}
\end{sub}
\vskip3mm

\begin{sub}\label{kern}{\bf Lemma:}\quad Let $X$ be an algebraic scheme over
  $k$ and $Y\subset X$ a closed subscheme defined by an ideal sheaf
  $\ki\subset\ko_X$ with $\ki^2=0$. Let $G$ be a linear algebraic group,
  let $\ko_X(G)$ denote the sheaf of morphisms from open sets of $X$
  to $G$, and let $L(G)$ denote the Lie algebra of $G$. Then there is a short 
exact sequence of sheaves of groups
\[
0\to L(G)\otimes_k\ki\to \ko_X(G)\to \ko_Y(G)\to 1.
\]

\begin{proof}
Since $G$ is smooth, $\ko_X(G)\rightarrow \ko_Y(G)$ is an epimorphism. In 
order to identify its kernel, we may assume that $G$ is a closed subgroup of 
${\text{GL}}_r(k)$ for some $r$. The group
${\text{GL}}_r(k)$ is an open subset of the 
affine space ${\text{Mat}}_r(k)$ of $r\times r$ matrices. 
Now, one has an exact sequence 
\[
0\to \Mat_r(k)\otimes_k\ki \overset{\varepsilon}{\lra} 
\ko_X({\text{GL}}_r(k))\to \ko_Y({\text{GL}}_r(k))\to 1,
\]
in which $\varepsilon$ is defined by $A\otimes f\mapsto e+Af$ as truncated
exponential, with $e$ denoting the unit $r\times r$ matrix. 
Let $I_G\subset k[(t_{ij})_{1\leq i,j\leq r}]$ be the ideal of polynomials
vanishing on $G$.
Then, for an element $\gamma\in\Mat_r(k)\otimes_k\ki (U)$,
where $U$ is open affine in $X$, $\varepsilon(\gamma)$
belongs to $\ko_X(G)(U)$ if and only if $F(\varepsilon(\gamma))=0$, 
for every polynomial $F\in I_G$. 
One may write 
$\gamma = A_1\otimes f_1+ \dots + A_m\otimes f_m$ with $A_1,\dots ,A_m\in 
\Mat_r(k)$ and with $f_1,\dots ,f_m\in \ki(U)$ linearly independent 
over $k$. 
Now, the Taylor expansion of any $F\in k[t_{ij}] $ at the identity 
$e\in \Mat_r(k),$ which reads as
\[
F(e+A\cdot f)=\sum\limits_{i,j}\frac{\partial F}{\partial t_{ij}}(e)\cdot 
a_{ij}\cdot f=(d_eF)(A)\cdot f,
\]
yields the formula
$F(\varepsilon(\gamma))=(\text{d}_eF)(A_1)\cdot f_1+\cdots +(\text{d}_eF)(A_m)
\cdot f_m,$
since $F(e)=0$ and $\ki(U)^2=0$.

If $\varepsilon(\gamma)\in \ko_X(G)(U)$, then 
it follows that $(\text{d}_eF)(A_\mu)=0$, $\mu=1,\dots ,m$, for any 
$F\in I_G.$ But the intersection of the 
kernels of the differentials $\text{d}_eF:{\text{Mat}}_r(k)\rightarrow k$ 
for all the $F\in I_G$ is exactly the tangent space ${\text{T}}_eG=L(G)$. 
Consequently, the kernel of $\ko_X(G)(U)\rightarrow \ko_Y(G)(U)$ is 
$L(G){\otimes}_k\ki(U)$. We have thus established the exact diagram 
\[
\xymatrix{0 \ar[r] & L(G)\otimes_k\ki\ar@{_{(}->}[d]\ar^{\varepsilon_G}[r] &
  \ko_X(G)\ar[r]\ar@{_{(}->}[d] & \ko_Y(G)\ar[r]\ar@{_{(}->}[d] & 1\\
0 \ar[r] & L(\GL_r)\otimes_k\ki\ar^{\varepsilon}[r] &
  \ko_X(\GL_r)\ar[r] & \ko_Y(\GL_r)\ar[r] & 1}
\]
with $\varepsilon_G$ induced by $\varepsilon$.  
\end{proof}
\end{sub}
\vskip3mm

\begin{sub}\label{equi}{\bf Remark:}\quad\rm
The action of $\ko_X(G)$ on itself deduced from the action of $G$ on 
itself by inner automorphisms induces, via the exact sequence from Lemma 
\ref{kern}, an action of $\ko_X(G)$ on $L(G){\otimes}_k\ki$. On the other 
hand, the action of $\ko_X(G)$ on $L(G){\otimes}_k\ko_X$ $($identified with  
the sheaf of morphisms from open sets of $X$ to the vector space $L(G)$$)$ 
deduced from the adjoint action of $G$ on $L(G)$ induces, via the exact 
sequence $:$
\[
0\to L(G){\otimes}_k\ki\to L(G){\otimes}_k\ko_X\to L(G){\otimes}_k\ko_Y\to 
0,
\]
an action of $\ko_X(G)$ on $L(G){\otimes}_k\ki$.
These two actions of $\ko_X(G)$ on $L(G){\otimes}_k\ki$ coincide since they 
obviously coincide in the case $G = {GL}_r$.
\end{sub}
\vskip5mm


 
\begin{sub}\label{restr}{\bf Lemma:}\quad
Under the assumptions of Lemma \ref{kern}, let $\mathcal F$ be a 
principal 
$G$-bundle over $X$ and let $\ke=\kf\vert Y$. Then there 
exists a canonical map 
$${\fam0 H}^1(Y,\ad(\ke)\otimes_{\ko_Y}\ki)\overset{\alpha}{\lra}
{\fam0 H}^1(X,\ko_X(G))$$ 
sending $0$ to the isomorphism class of $\kf$, and whose image
is the set of isomorphism classes of 
principal $G$-bundles $\kf^{\prime}$ over $X$ such that 
$\kf^{\prime}\vert Y\simeq\ke$. 

\begin{proof}
$\kf$ corresponds to an element $c\in \text{H}^1(X,\ko_X(G))$. 
If one uses the adjoint action of $\ko_X(G)$ on 
$L(G)\otimes_k\ko_X$, then the corresponding 
twisted sheaf $(L(G)\otimes_k\ko_X)^c$ is 
exactly $\text{ad}({\mathcal F})$. The conclusion follows now from 
Lemma \ref{kern} and Lemma \ref{idim}, taking into account that, according 
to the above Remark \ref{equi}, one has an exact sequence :
$$0\rightarrow (L(G)\otimes_k\ki)^c\rightarrow 
(L(G)\otimes_k\ko_X)^c\rightarrow 
(L(G)\otimes_k\ko_Y)^c\rightarrow 0,$$
hence : 
$(L(G)\otimes_k\ki)^c\simeq 
\text{Ker}(\text{ad}(\kf)\rightarrow \text{ad}(\kf)\vert Y) 
\simeq \text{ad}(\kf)\otimes_{{\ko}_X}\ki\simeq 
\text{ad}(\ke)\otimes_{\ko_Y}\ki$. 
\end{proof}

\rm Notice, for further use, that, by construction, \textit{the map} 
$\alpha $   
\textit{in the statement of the previous lemma is functorial in} 
$(X,Y,{\kf})$. 
\end{sub}
\vskip10mm

\section{Proof of Theorem 1}

First, let us recall a result which is implicit in Kempf's paper 
\cite{k}. For an explicit proof see \cite{ct}.

\begin{sub}\label{kempf}{\bf Lemma:}\quad
Let $E$ be an algebraic vector bundle on $\P_n$, 
$n\geq 2$, $H\subset \P_n$ a hyperplane, $x\in \P_n
\setminus H$ and $p:\P_n\smallsetminus \{x\}\rightarrow H$ the central 
projection. If $E$ and $p^{\ast}(E|H)$ are isomorphic, as vector bundles, 
over each infinitesimal neighborhood of $H$ in $\P_n$, then $E$ 
splits into a direct sum of line bundles.
\end{sub}
\vskip3mm

In characteristic 0 one can generalise the above lemma to principal 
bundles : 

\begin{sub}\label{kempfpb}{\bf Lemma:}\quad 
Assume that ${\fam0 char}(k)=0$ and let $G$ be a linear algebraic 
group over $k$. 
Let $\ke$ be a principal $G$-bundle on $\P=\P_n, n\geq 2$,
and let $H$ and $p$ be as in the previous lemma.
If $\ke$ and $p^\ast(\ke|H)$ are isomorphic as principal
$G$-bundles over each infinitesimal neighborhood of $H$ in $\P_n$, then
$\ke$ is split.

\begin{proof} Let $c\in {\text{H}}^1(\P, \ko_{\P}(G))$ be the class of
$\ke$. Let $R_uG$ be the unipotent radical of $G$, $Q=G/R_uG$ the reductive 
quotient and $\rho : G\rightarrow Q$ the canonical surjection. We will show 
that ${\text{H}}^1(\P,\ko_{\P}(R_uG)^c)=\{e\}$ (see (II) below) and that 
$(L(Q){\otimes}_k\ko_{\P})^c$ is a direct sum of line bundles (see (I) below). 

Now, according to a result of G.D. Mostow \cite{mo} (which is valid only in 
characteristic 0, see \cite{bs} or \cite{bt}) there exists a 
\textit{Levi subgroup} $\Lambda $ of $G$, i.e., a closed subgroup such that, 
denoting by $u$ the inclusion $\Lambda \hookrightarrow G$,
the composition $\rho \circ u : 
\Lambda \rightarrow Q$ is an isomorphism. Since 
${\text{H}}^1(\P,\ko_{\P}(R_uG)^c)=\{e\}$, Lemma \ref{idim} implies that 
${\text{H}}^1(\rho )^{-1}({\text{H}}^1(\rho )(c))=\{c\}$. In particular, if 
$c_{\Lambda}\in {\text{H}}^1(\P,\ko_{\P}(\Lambda ))$ is defined by 
${\text{H}}^1(\rho \circ u)(c_{\Lambda})={\text{H}}^1(\rho )(c)$ then 
$c={\text{H}}^1(u)(c_{\Lambda})$, i.e., $\ke$ admits a reduction of structure 
group to $\Lambda $. Let $\ke_{\Lambda}$ be the principal $\Lambda $-bundle 
defined by $c_{\Lambda}$. Then 
$$\text{ad}(\ke_{\Lambda}):=(L(\Lambda ){\otimes}_k\ko_{\P})^{c_{\Lambda}}
\simeq (L(Q){\otimes}_k\ko_{\P})^c$$ 
is a direct sum of line bundles. Since $\Lambda \simeq Q$ is reductive and 
$\text{char}(k)=0$, \cite{BGH}, Theorem 4.3., implies that $\ke_{\Lambda}$ 
admits a reduction of structure group to a maximal torus $T$ of $\Lambda $.

This proves the lemma, modulo the two technical facts (I) and (II) quoted 
above. 

(I) For any closed normal connected subgroup $N$ of $G$ there is 
the induced exact sequence of Lie algebras :
\[
0\to L(N)\to L(G)\to L(G/N)\to 0
\]
and the associated exact sequence of locally free sheaves
\[
0\to (L(N){\otimes}_k\ko_{\P})^c\to (L(G){\otimes}_k\ko_{\P})^c\to 
(L(G/N){\otimes}_k\ko_{\P})^c\to 0.
\]

One sees easily that each of the three bundles occurring in the last exact
sequence satisfies the hypothesis of Lemma \ref{kempf}, hence is a direct sum
of line bundles. Moreover, since $n\geq 2$, this exact sequence splits, so
that $(L(N){\otimes}_k\ko_{\P})^c$ is a direct summand of 
$(L(G){\otimes}\ko_{\P})^c$.
\newpage

(II) We show now that ${\text{H}}^1(\P, \ko_{\P}(R_uG)^c)=\{e\}$. 

To prove that we consider the central series
\[
R_uG=C^0\supset C^1\supset\ldots\supset C^n=\{e\}
\]
of $R_uG$. Each of the groups $C^i$ is a closed connected normal subgroup of
$G$ and the quotients $C^i/C^{i+1}$ are abelian and unipotent. This implies
that the exponential map
\[
L(C^i/C^{i+1})\rightarrow C^i/C^{i+1}
\]
is an isomorphism of algebraic groups.
Using again the twisting by $c$, which is induced by the inner automorphisms 
of $G$, we obtain the exact sequences
\[
0\to (L(C^{i+1}){\otimes}_k\ko_{\P})^c\to (L(C^i){\otimes}_k\ko_{\P})^c\to
(L(C^i/C^{i+1}){\otimes}_k\ko_{\P})^c\to 0
\]
and, according to (I), $(L(C^{i+1}){\otimes}_k\ko_{\P})^c$ is a direct summand
to $(L(C^i){\otimes}_k\ko_{\P})^c$ for $i\geq 0$.

It follows that also $(L(C^i/C^{i+1}){\otimes}_k\ko_{\P})^c$ is a direct sum 
of line bundles. Since $n\geq 2$, 
${\text{H}}^1(\P, (L(C^i/C^{i+1}){\otimes}_k\ko_{\P})^c)=0$
for $i\geq 0$, and then also ${\text{H}}^1(\P,\ko_{\P}(C^i/C^{i+1})^c)=\{e\}$. 
This proves that
${\text{H}}^1(\P,\ko_{\P}(R_uG)^c)=\{e\}$.
\end{proof}
\end{sub}
\vskip8mm

We are able, now, to prove Theorem 1.
Using the notation from the preparations preceding the 
proof of the Theorem in \cite{ct}, suppose that there exists a principal 
$G$-bundle $\mathcal F$ over ${\mathbb P}_{n+m}$ such that 
$\kf\vert L\simeq \mathcal E$. We shall construct a homogeneous 
ideal $J\subset R$, generated by 
${\Sigma}_{i>0}\text{dimH}^1(\text{ad}(\mathcal E)(-i))$ homogeneous elements 
such that, for any $i\geq 0$, 
$$\kf\vert L_i\cap X\simeq 
{\pi}^{\ast}\ke\vert L_i\cap X,$$ 
where $X$ is the closed subscheme 
of ${\mathbb P}_{n+m}$ defined by the ideal $JS$ and $\pi : 
{\mathbb P}_{n+m}\setminus L^{\prime}\rightarrow L$ the central projection.
The inequality imposed on $m$ implies that there exists $p\in L^{\prime}
\simeq {\mathbb P}_{m-1}$ such that the polynomials from $J$ vanish in $p$. 
The linear span $P$ of $p$ and $L$ is contained in $X$, hence 
$\kf\vert L_i\cap P\simeq {\pi}^{\ast}\ke\vert L_i\cap P$, 
for any $i\geq 0$. 
 
Recall that $P\simeq {\P}_{n+1}$ and that the schemes 
$L_i\cap P$ are the infinitesimal neighborhoods in $P$ of the hyperplane 
$L$ of $P$. 
Therefore, if $\text{char}(k)=0$, Lemma \ref{kempfpb}
implies that $\kf\vert P$ is split and so is $\ke\simeq \kf\vert L$.

In the case of a reductive linear algebraic group in arbitrary 
characteristic, we know that also
$\text{ad}(\kf)\vert L_i\cap P
\simeq {\pi}^{\ast}\text{ad}(\ke)\vert L_i\cap P$ (as vector bundles 
on $L_i\cap P$), and then Lemma \ref{kempf} implies that $\text{ad}(\ke)$
splits as a direct sum of line bundles. From Proposition 1 one deduces that 
$\ke$ is split in this case, too.


Finally, $J$ is constructed, as in the proof of the Theorem in 
\cite{ct}, by a standard technique borrowed from infinitesimal deformation 
theory, using Lemma \ref{restr} above. Explicitly:

Suppose that $J\subset R$ has already been constructed such that
$\kf|L_i\cap X\simeq \pi^\ast\ke|L_i\cap X$. We enlarge $J$ in
degree $\geq i+1$ to an ideal $J'$ as to obtain also 
$\kf|L_{i+1}\cap X'\simeq\pi^\ast\ke| L_{i+1}\cap X',$ 
where $X^{\prime}$ is the closed subscheme of ${\mathbb P}_{n+m}$ defined by 
the ideal $J^{\prime}S$.

To do so we put $X_j=L_j\cap X.$ 
Using the notation of \cite{ct}, the ideal sheaf $\ki_{X_i}$ of $X_i$ in
$X_{i+1}$ is isomorphic to $\ko_L(-i-1)\otimes R_{i+1}/J_{i+1}$ and satisfies
$\ki^2_{X_i}=0$. By Lemma \ref{restr} there is a canonical map
\[
{\text{H}}^1(L, \ad (\ke)(-i-1)\otimes R_{i+1}/J_{i+1})\overset{\alpha}{\lra}
{\text{H}}^1(X_{i+1}, \ko_{X_{i+1}}(G)) 
\]
such that $\alpha(0)=[\pi^\ast\ke|X_{i+1}]$ and 
the image of $\alpha$ is the set of all classes $[\kf']$ of
principal bundles $\kf'$ on $X_{i+1}$ such that
$\kf'|X_i\simeq\pi^\ast\ke|X_i$. By assumption the class of
$\kf|X_{i+1}$ belongs to this set, hence  
$[\kf|X_{i+1}]=\alpha(\xi)$ for some $\xi\in {\text{H}}^1(L, \ad(\ke)(-i-1))
\otimes  R_{i+1}/J_{i+1}$. Let $\xi_1, \ldots, \xi_s$ be a basis of 
${\text{H}}^1(L,\ad(\ke)(-i-1))$. Then
\[
\xi=\xi_1\otimes \bar{f_1}+\cdots+\xi_s\otimes \bar{f_s}
\]
with unique residue classes $\bar{f_\nu}\in R_{i+1}/J_{i+1},\; f_\nu\in
R_{i+1}$. Let
\[
J':=J+Rf_1+\cdots+Rf_s
\]
and let $X'\subset X$ be the variety of $J'S\supset JS$.

Then
$X'_i=L_i\cap X'=L_i\cap X=X_i$
and $X_i\subset X'_{i+1}\subset X_{i+1}$. 

According to the functoriality of the maps $\alpha$ in Lemma 
\ref{restr}, there is a commutative diagram
\[
\xymatrix{{\text{H}}^1(L, \ad(\ke)(-i-1))\otimes
  R_{i+1}/J_{i+1}\ar^-{\alpha}[r]\ar^{\rho'}[d] & {\text{H}}^1(X_{i+1},
  \ko_{X_{i+1}}(G))\ar^\rho[d] \\
{\text{H}}^1(L, \ad(\ke)(-i-1))\otimes  
R_{i+1}/J'_{i+1} \ar^-{\alpha^{\prime}}[r]            & {\text{H}}^1(X'_{i+1},
\ko_{X'_{i+1}}(G))
},
\]
where $\rho'$ and $\rho$ denote the natural quotient maps.
By definition of $\alpha^{\prime}$ in Lemma \ref{restr}, 
$[\pi^\ast\ke|X'_{i+1}]=\alpha^{\prime}(0)$. 
Since ${\rho}^{\prime}(\xi )=0$, it follows that
\[
[\kf|X'_{i+1}]=[\pi^\ast\ke|X'_{i+1}].
\]
This completes the inductive construction of $J$ and the proof of Theorem 1.
\vskip10mm


\section{Proof of Theorem 2}
 We use a trick of Mohan Kumar \cite{mk}, to show that, under the 
hypothesis of the theorem, 
$\text{ad}(\mathcal E)$ can be extended to a vector bundle on 
${\mathbb P}_{n+1}$.

\par Embed ${\mathbb P}_n$ as the hyperplane $H$ of ${\mathbb P}_{n+1}=:P$ 
of equation $X_{n+1}=0$ and let $H_i$ denote its $i$th infinitesimal 
neighbourhood, of equation $X_{n+1}^{i+1}=0$. Let $x\in P\setminus H$ and let 
${\pi}_x : P\setminus \{x\}\rightarrow H$ be the projection. Using Lemma 
\ref{restr} and the vanishing conditions in the hypothesis 
one shows, by induction on $i\geq 0$, that if $\mathcal F$ is a principal 
$G$-bundle over $H_i$ such that ${\mathcal F}\vert H\simeq \mathcal E$ then 
${\mathcal F}\simeq {\pi}_x^{\ast}{\mathcal E}\vert H_i$. In particular, 
if $y\in P\setminus H$ is another point and ${\pi}_y :P\setminus \{y\} 
\rightarrow H$ the corresponding projection, then 
${\pi}_y^{\ast}{\mathcal E}\vert H_i\simeq {\pi}_x^{\ast}{\mathcal E}\vert 
H_i$, $\forall i\geq 0$. This implies that 
${\pi}_y^{\ast}\text{ad}(\mathcal E)\vert H_i\simeq 
{\pi}_x^{\ast}\text{ad}(\mathcal E)\vert H_i$, $\forall i\geq 0$.

\par Both ${\pi}_x^{\ast}\text{ad}(\ke)$ and 
${\pi}_y^{\ast}\text{ad}(\ke)$ can be extended to {\it reflexive sheaves} 
${\mathcal A}_x$ and ${\mathcal A}_y$ on $P$. The sheaf
${\mathcal Hom}_{{\mathcal O}_P}({\mathcal A}_x,{\mathcal A}_y)$ is 
reflexive, hence, for $j=0,\  1$, 
$\text{H}^j({\mathcal Hom}_{{\mathcal O}_P}({\mathcal A}_x,{\mathcal A}_y)
(-i-1))=0$ for $i>>0$. It follows that 
$$\text{Hom}_{{\mathcal O}_P}({\mathcal A}_x,{\mathcal A}_y)
\overset{\sim}{\longrightarrow} 
\text{Hom}_{{\mathcal O}_{H_i}}({\mathcal A}_x\vert H_i,
{\mathcal A}_y\vert H_i)$$
for $i>>0$. For $i>>0$, any isomorphism ${\mathcal A}_x\vert H_i
\overset{\sim}{\rightarrow} {\mathcal A}_y\vert H_i$ can be lifted to a 
morphism ${\mathcal A}_x\rightarrow {\mathcal A}_y$ which must be an 
isomorphism on a (Zariski) open neighbourhood $U$ of $H$ in $P$, that is, 
${\pi}_x^{\ast}\text{ad}(\ke)$ and ${\pi}_y^{\ast}\text{ad}(\ke)$ 
are isomorphic over $U$.  
But $P\setminus U$ must be 0-dimensional, hence it has codimension $\geq 2$. 
It follows that ${\pi}_x^{\ast}\text{ad}(\ke)$ and 
${\pi}_y^{\ast}\text{ad}(\ke)$ 
are isomorphic over $P\setminus \{x,y\}$, hence they 
can be glued and one gets a vector bundle $\tilde A$ on $P$ extending 
$\text{ad}(\ke)$. 

\par Since ${\tilde A}\vert H_i\simeq 
{\pi}_x^{\ast}\text{ad}(\ke)\vert H_i$, $\forall i\geq 0$,  
Lemma \ref{kempf}  implies that 
$\tilde A$ splits, hence $\text{ad}(\ke)\simeq {\tilde A}\vert H$ splits.
\vskip3mm

\section{Proof of Proposition 1.1}\label{pzartop}
A theorem of Chevalley says that there is
a short exact sequence of groups
\begin{equation}\label{e1pn}
1\longrightarrow H  \longrightarrow G \longrightarrow A
\longrightarrow 1\, ,
\end{equation}
where $A$ is an abelian variety over $k$
and $H$ is an affine algebraic
group over $k$; for a modern proof see \cite{bc}. We will show that any 
algebraic principal $A$-bundle over ${\mathbb P}_n$ is trivial.

Let $E_A$ be a principal $A$-bundle over ${\mathbb P}_n$.
To prove that $E_A$ is trivial, it suffices to 
show that $E_A$ admits a
section over the generic point. Indeed, if $s$ is a section of $E_A$
over a Zariski open subset of ${\mathbb P}_n$, then $s$ extends to a 
section of the pullback of $E_A$
over some blow-up of ${\mathbb P}_n$. Since an abelian variety does not 
have any rational curves, the section over the blow-up of ${\mathbb 
P}_n$ descends to a section of $E_A$ over ${\mathbb P}_n$.

There is a separable extension $K'$ of the function field $K$
of ${\mathbb P}_n$ over which $E_A$ has a rational point. Hence
$E_A$ over $K'$ is trivial. There is an inflation
homomorphism $\text{H}^1(K', A) \longrightarrow \text{H}^1(K, A)$ whose
composition with the natural homomorphism $\text{H}^1(K, A) \longrightarrow
\text{H}^1(K', A)$ is multiplication by $d$, where $d$ is the degree
of the field-extension. So the class in $\text{H}^1(K, A)$ given by
$E_A$ is torsion. We noted earlier that a principal $A$-bundle
over ${\mathbb P}_n$ is trivial if its restriction to $K$ is trivial.
Therefore, the class in $\text{H}^1(K, A)$ given by
$E_A$ being torsion it follows that the class in 
$\text{H}^1({\mathbb P}_n, A)$
given by $E_A$ is torsion. Consequently, the principal $A$-bundle
$E_A$ over ${\mathbb P}_n$ admits a reduction of
structure group to a finite group-scheme. Since the fundamental
group-scheme of ${\mathbb P}_n$ is trivial \cite[p. 93, Corollary]{No},
it follows that any principal bundle over ${\mathbb P}_n$ with a finite
group-scheme as the structure group is trivial. Hence $E_A$ is
trivial.

Since any principal $A$-bundle over ${\mathbb P}_n$ is trivial,
using (\ref{e1pn}) it follows that any principal $G$-bundle
over ${\mathbb P}_n$ admits a reduction of structure group to
the subgroup $H$. Therefore, to prove the Proposition
it suffices to show that any principal $H$-bundle over
${\mathbb P}_n$ is Zariski locally trivial.

Let $E_H$ be a principal $H$-bundle over ${\mathbb P}_n$.
Let $H_0\, \subset\, H$ be the connected component of $H$ containing
the identity element. Let $E_{H/H_0} = E_H\times_H (H/H_0)$ be
the principal $(H/H_0)$-bundle over ${\mathbb P}_n$ obtained
by extending the structure group of $E_H$ using the
quotient map
$H\longrightarrow H/H_0$. Since ${\mathbb P}_n$ is simply
connected, it follows immediately that $E_{H/H_0}$ is a trivial
principal $(H/H_0)$-bundle. Therefore, $E_H$ admits a reduction
of structure group to $H_0$.

Let $E_{H_0}$ be a principal $H_0$-bundle over ${\mathbb P}_n$.
To prove the Proposition it is enough to show that
$E_{H_0}$ is Zariski locally trivial.

We will prove that $H_0$ is acceptable in the sense of
\cite[p. 188, Definition]{Ra}. But before that we will
show that $E_{H_0}$ is Zariski locally trivial assuming that
$H_0$ is acceptable.

So assume that $H_0$ is acceptable. Since $k$ is algebraically
closed, any principal $H_0$-bundle over $\text{Spec}\, k$
is trivial. Hence using \cite[p. 189, Theorem A]{Ra} it
follows that the restriction of $E_{H_0}$ to some open subscheme
of any affine chart of $\P_n$ is trivial.


It follows from \cite[hypothesis (1), p. 97]{CO} or 
\cite[p. 110, Theorem 3.2]{CO}, and the assumption that 
$k$ is algebraically closed, that the
principal $H_0$-bundle $E_{H_0}$ is Zariski locally trivial
under the assumption that $H_0$ is acceptable.

To prove that $H_0$ is acceptable, let $R_uH_0$ be the
unipotent radical of $H_0$. So we have a short exact sequence
of groups
$$
1\longrightarrow R_uH_0 \longrightarrow H_0
\longrightarrow Q_0\longrightarrow 1\, ,
$$
where $Q_0$ is reductive. Note that $Q_0$ is connected
as $H_0$ is so. From \cite[p. 137, Theorem 1.1]{RR}
we know that $Q_0$ is acceptable. Hence it suffices to
show that $R_uH_0$ is acceptable.

The unipotent group $R_uH_0$ has a filtration of normal
subgroups
$$
e= U_0 \subset U_1\subset \cdots \subset U_i 
\subset \cdots \subset U_{d-1} \subset U_d= R_uH_0\, ,
$$
where $d= \dim R_uH_0$, and $U_i/U_{i-1}$ is the additive
group ${\mathbb G}_a$ for each $i\in [1 ,d]$ (see \cite[p. 123,
Theorem 19.3]{hu}). Therefore, the group $R_uH_0$
is acceptable if ${\mathbb G}_a$ is acceptable. But
$$
\text{H}^1_{et}({\mathbb A}^1, {\mathbb G}_a) = 
\text{H}^1_{et}({\mathbb A}^1,{\mathcal O}_{{\mathbb A}^1}) =0\, .
$$
Hence ${\mathbb G}_a$ is acceptable. This completes the proof
of Proposition \ref{zartop}.

\section{Proof of Proposition 1}\label{prop1}

The aim in this section is to prove Proposition 1 for algebraically
closed fields $k$ of arbitrary characteristics.
(for $k = \C$ this follows from \cite{BGH})

\vskip3mm
{\bf Proposition:}\quad
\textit{Let $G$ be reductive linear algebraic group
defined over $k$. Let $\ke$ be a principal $G$-bundle over ${\P}_n$
and $\ad(\ke)$ its adjoint 
bundle. If $\ad(\ke)$ splits as 
a direct sum of line bundles, then $\ke$ is split.}
\vskip3mm

\begin{proof}
We recall that if the characteristic of the base field
$k$ is positive, then a principal bundle $\mathcal F$ over a
smooth variety $X$ defined over $k$
is called \textit{strongly semistable}
if the pull back $(F^n_X)^*\mathcal F$ over $X$ is semistable
for all $n\geq 1$, where $F^n_X$ is the n-fold composition
of the Frobenius morphism $F_X : X\rightarrow X$. For
our convenience, when the characteristic of $k$ is zero, by
a strongly semistable principal bundle we will mean a
semistable bundle.

Since the tangent bundle $T{\P}_n$ is semistable of positive
degree, any semistable vector bundle over ${\P}_n$ is
strongly semistable \cite[p. 316, Theorem 2.1(1)]{MR}.

Let $\ke$ be a principal $G$-bundle over $T{\P}_n$, where $G$
is a reductive linear algebraic group defined over $k$. Then
$\ke$ admits a unique Harder-Narasimhan reduction; see \cite{BH}. 
In general some conditions are needed for the uniqueness part
of the Harder-Narasimhan reduction. (See \cite[p. 208,
Proposition 3.1]{BH} for the existence of Harder-Narasimhan reduction,
and \cite[p. 221, Corollary 6.11]{BH} for the uniqueness; the
fact that any any semistable vector bundle over ${\P}_n$ is
strongly semistable ensures Proposition 6.9 in  \cite[p. 219]{BH}
remains valid without the assumption on the height.)

Now we assume that the adjoint vector bundle $\ad(\ke)$ is
a direct sum of line bundles. This immediately implies that the
Harder-Narasimhan filtration of $\ad(\ke)$ is a filtration
of subbundles of $\ad(\ke)$ (in general it is only a filtration
of subsheaves with each successive quotient being torsionfree).
Therefore, the Harder-Narasimhan reduction of $\ke$ is defined
over entire ${\P}_n$.

Let
\begin{equation}\label{red.}
{\ke}_P\subset \ke
\end{equation}
be a principal $P$-bundle giving the Harder-Narasimhan reduction
of $\ke$ over ${\P}_n$; here $P\subset G$ is a parabolic subgroup.

Any principal $G$-bundle over ${\P}_1$ is split
(see \cite{Ha}). Therefore, the proposition is proved
for $n =1$. Henceforth, we will assume that $n \geq 2$.

Consider the short exact sequence of groups
\begin{equation}\label{sh.ex}
1\rightarrow R_u(P) \rightarrow P \rightarrow Q(P)\rightarrow 1\, ,
\end{equation}
where $R_u(P)$ is the unipotent radical of $P$, and $Q(P)$
is the Levi quotient of of $P$. This short exact sequence is
right split. Fix a subgroup of $P$ that projects isomorphically
to $Q(P)$. This subgroup will be denoted by $\Lambda (P)$.
We will show that ${\ke}_P$ admits a reduction of structure
group to the subgroup $\Lambda (P)$ of $P$.

To prove this first note that giving a reduction of structure
group of ${\ke}_P$ to $\Lambda (P)$ is equivalent to giving a section
of the fibre bundle ${\ke}_P/\Lambda (P)$ over ${\P}_n$.
Let ${\ke}_P(R_u(P))$ be the group-scheme over ${\P}_n$
associated to ${\ke}_P$ for the adjoint action of $P$ on the
normal subgroup $R_u(P)$ in \eqref{sh.ex}. The fibre bundle
${\ke}_P/\Lambda (P)$ is a torsor for ${\ke}_P(R_u(P))$. In other words,
the fibres of ${\ke}_P(R_u(P))$ have a natural free transitive
action on the fibres of ${\ke}_P/\Lambda (P)$. Torsors for
${\ke}_P(R_u(P))$ are parametrised by $\text{H}^1({\P}_n,
{\ke}_P(R_u(P)))$. Therefore, to prove that ${\ke}_P$
admits a reduction of structure
group to the subgroup $\Lambda (P)$ it suffices to show that
\begin{equation}\label{eq00}
\text{H}^1({\P}_n, {\ke}_P(R_u(P)))= 0\, .
\end{equation}

Consider the upper central series $\{{\mathcal G}_i\}_{i\geq 0}$
for $R_u(P)$. So ${\mathcal G}_0 = R_u(P)$ and
${\mathcal G}_{i+1} = [R_u(P), {\mathcal G}_i]$ for all
$i \geq 0$. This central series is preserved by the adjoint
action of $P$, and each successive quotient is an abelian
unipotent group. For an abelian unipotent group, the exponential
map from its Lie algebra is well-defined, and it is an isomorphism.
Also, for any line bundle $\xi$ on ${\P}_n$ we have
$\text{H}^1({\P}_n, \xi) =0$ (recall that $n>1$). Therefore it follows that
\eqref{eq00} holds.

Let
\begin{equation}\label{lp}
{\ke}_{\Lambda (P)}\subset {\ke}_P
\end{equation}
be a reduction of structure group of ${\ke}_P$ to $\Lambda (P)$.
Let ${\ke}'_{Q(P)}$ be the principal $Q(P)$-bundle obtained
by extending the structure group of ${\ke}_P$ using the
projection in \eqref{sh.ex}. We note that the principal
$Q(P)$-bundle ${\ke}'_{Q(P)}$ is
canonically identified with ${\ke}_{\Lambda (P)}$.

Let
\[
F_1 \subset \cdots \subset F_{m-1} \subset F_m = \ad({\ke})
\]
be the Harder-Narasimhan filtration of $\ad({\ke})$. From the
construction of the Harder-Narasimhan reduction of ${\ke}$
we know that $m=2m_0+1$, and
$\text{degree}(F_{(m+1)/2}/F_{(m-1)/2})=0$. Furthermore,
the subbundle $\text{ad}({\ke}_{P})\subset \ad({\ke})$ coincides
with $F_{(m+1)/2}$, and $\text{ad}({\ke}'_{Q(P)})$ coincides
with the quotient $F_{(m+1)/2}/F_{(m-1)/2}$. In particular,
$\text{ad}({\ke}'_{Q(P)})$ is a semistable vector bundle.
We noted earlier that any semistable vector bundle on ${\P}_n$ is
strongly semistable. Therefore, $\text{ad}({\ke}'_{Q(P)})$
is strongly semistable.

Since ${\ke}'_{Q(P)}$ is identified with ${\ke}_{\Lambda (P)}$, we now
conclude that the adjoint vector bundle
$\text{ad}({\ke}_{\Lambda (P)})$ is strongly semistable.

We recall that $\ad(\ke)$
is a direct sum of line bundles. From the above remark
that the subbundle
$\text{ad}({\ke}_{P})$ coincides with $F_{(m+1)/2}$ it follows
immediately that $\ad({\ke}_P)$
is also a direct sum of line bundles. Since the adjoint bundle
$\ad({\ke}_{\Lambda (P)})$ is a direct summand of $\ad({\ke}_P)$, using
the Atiyah-Krull-Schmidt theorem (see
\cite[p. 315, Theorem 3]{At}) we conclude that
$\ad({\ke}_{\Lambda (P)})$ is also a direct sum of line bundle.

Since $\Lambda (P)$ is reductive, the adjoint group
\begin{equation}\label{ad.lp}
H := \Lambda (P)/Z(\Lambda (P))
\end{equation}
$\Lambda (P)$ does not admit any nontrivial character; here
$Z(\Lambda (P))$ is the center of $\Lambda (P)$.
Hence $\det \ad({\ke}_{\Lambda (P)})=\bigwedge^{\text{top}}
\ad({\ke}_{\Lambda (P)})$ is a trivial line bundle. On the other
hand, we already proved that $\ad({\ke}_{\Lambda (P)})$ is semistable
and it splits into a direct sum of line bundles. Combining
these it follows that
$\ad({\ke}_{\Lambda (P)})$ is a trivial vector bundle.

Let ${\mathfrak l}({\mathfrak p})$
denote the Lie algebra of the group $\Lambda (P)$. Consider the
adjoint action of $\Lambda (P)$ on ${\mathfrak l}({\mathfrak p})$.
It gives a homomorphism to the linear group
\begin{equation}\label{rho}
\rho : \Lambda (P)\rightarrow \text{GL}({\mathfrak l}({\mathfrak p}))\, .
\end{equation}
Let ${\ke}_{\text{GL}({\mathfrak l}({\mathfrak p}))}$ denote the
principal $\text{GL}({\mathfrak l}({\mathfrak p}))$-bundle
over ${\P}_n$ obtained by extending the structure group of
${\ke}_{\Lambda (P)}$ using the homomorphism $\rho$ in \eqref{rho}.
We noted earlier that $\ad({\ke}_{\Lambda (P)})$ is a trivial vector
bundle. Therefore, ${\ke}_{\text{GL}({\mathfrak l}({\mathfrak p}))}$
is a trivial principal bundle.

Consider the quotient $H$ of $\Lambda (P)$ defined in \eqref{ad.lp}.
Let ${\ke}_{H}$ be the principal $H$-bundle over ${\P}_n$
obtained by extending the structure group of ${\ke}_{\Lambda (P)}$
using the quotient map. The homomorphism $\rho$
in \eqref{rho} factors through $H$. Therefore,
${\ke}_{\text{GL}({\mathfrak l}({\mathfrak p}))}$ is an extension
of structure group of ${\ke}_{H}$.

Since $\Lambda (P)$ is reductive, the homomorphism $\rho$
gives an embedding of $H$ into
$\text{GL}({\mathfrak l}({\mathfrak p}))$. We already noted
that ${\ke}_{\text{GL}({\mathfrak l}({\mathfrak p}))}$
is trivial. Therefore, the reduction 
\[
{\ke}_{H}\subset {\ke}_{\text{GL}({\mathfrak l}({\mathfrak p}))}
\]
is given by a morphism
\begin{equation}\label{f}
f :  {\P}_n\rightarrow \text{GL}({\mathfrak
l}({\mathfrak p}))/H\, .
\end{equation}

Since $H$ is a reductive subgroup of
$\text{GL}({\mathfrak l}({\mathfrak p}))$, the quotient space
$\text{GL}({\mathfrak
l}({\mathfrak p}))/H$ is an affine variety. Therefore, the
morphism $f$ in \eqref{f} is a constant one. This immediately
implies that the principal $H$-bundle ${\ke}_{H}$ is
trivial. From this it follows that the principal $\Lambda (P)$-bundle
${\ke}_{\Lambda (P)}$ admits a reduction of structure group to the
center $Z(\Lambda (P))$.

Since ${\ke}$ is an extension of structure group of ${\ke}_{\Lambda (P)}$,
we now conclude that ${\ke}$ admits a reduction of structure
group to a maximal torus of $G$.
\end{proof}


\begin{thebibliography}{888}
\bibitem{At} \textbf{M.F.~Atiyah}, On the Krull-Schmidt theorem with
application to sheaves, Bull. Soc. Math. Fr. 84 (1956), 307--317
 
\bibitem{b} \textbf{I.~Biswas}, A criterion for the splitting of a principal
bundle over a projective space, Arch. Math. 84 (2005), 38--45
 
\bibitem{BGH} \textbf{I.~Biswas, T.L. G\'omez, Y.I.~Holla},
Reduction of structure group of principal bundles over projective
manifolds with Picard number one, Int. Math. Res. Not. (2002),
889--903
 
\bibitem{BH} \textbf{I.~Biswas, Y.I.~Holla}, Harder-Narasimhan reduction
of a principal bundle, Nagoya Math. Journ. 174 (2004), 201--223
 
\bibitem{bs} \textbf{A.~Borel, J.-P.~Serre}, Th\' eor\` emes de
finitude en cohomologie galoisienne, Comm. Math. Helv. 39 (1964),
111--164
 
\bibitem{bt} \textbf{A.~Borel, J.~Tits}, Groupes reductifs,
Publ. Math. I.H.E.S. 27 (1965), 55--150
 
\bibitem{ct} \textbf{I.~Coand\u a, G.~Trautmann}, The splitting criterion of
Kempf and the Babylonian tower theorem, Comm. in Algebra 34 (2006),
2485--2488
 
\bibitem{CO} \textbf{J.-L.~Colliot-Th\' el\` ene, M.~Ojanguran},
Espaces principaux homog\` enes localement triviaux,
Publ. Math. I.H.E.S. 75 (1992), 97--122
 
\bibitem{bc} \textbf{B.~Conrad}, A modern proof of Chevalley's theorem on
  algebraic groups, Journ. Ramanujan Math. Soc. 17 (2002), 1--18
 
\bibitem{fr} \textbf{J.~Frenkel}, Cohomologie non ab\' elienne et espaces
fibr\' es, Bull. Soc. Math. France 85 (1957), 135--220
 
\bibitem{g} \textbf{A. Grothendieck}, Sur la classification des
fibr\' es holomorphes sur la sph\' ere de Riemann, Amer. Jour. Math. 79
(1957), 121--138

\bibitem{Ha} \textbf{G.~Harder}, Halbeinfache Gruppenschemata
\" uber vollst\" andigen Kurven, Invent. Math. 6 (1968), 107--149

\bibitem{hu} \textbf{J.E.~Humphreys}, Linear Algebraic Groups,
Springer Verlag, 1975
 
\bibitem{k} \textbf{G.R.~Kempf}, A criterion for the splitting of a vector
bundle, Forum Math. 2 (1990), 477--480
 
\bibitem{mk} \textbf{N.~Mohan~Kumar}, Vector bundles on projective spaces,
In : C. Musili (ed.) et al., Advances in algebra and geometry, Proceedings
of the international conference on algebra and geometry, Hyderabad, India,
December 7--12, 2001, pp. 185--188,
Hindustan Book Agency, New Delhi, 2003
 
\bibitem{MR} \textbf{V.B.~Mehta, A.~Ramanathan}, Homogeneous bundles
in characteristic $p$, In: Algebraic Geometry---open problems
(Ravello, 1982), pp. 315--320
Lect. Notes Math., 997, Springer, Berlin, 1983
 
\bibitem{mo} \textbf{G.D.~Mostow}, Fully reducible subgroups of algebraic
groups, Amer. Journ. Math. 78 (1956), 401--443
 
\bibitem{No} \textbf{M.V.~Nori}, The fundamental group scheme,
Proc. Ind. Acad. Sci. (Math. Sci.) 91 (1982), 73--122.
 
\bibitem{Ra} \textbf{M.S.~Raghunathan}, Principal bundles on affine
space. In: C. P. Ramanujam---a tribute, pp. 187--206,
Springer, Berlin-New York, 1978
 
\bibitem{RR} \textbf{M.S.~Raghunathan}, \textbf{A.~Ramanathan},
Principal bundles on the affine line, Proc. Ind. Acad. Sci. (Math. Sci.)
93 (1984), 137--145.



\end{thebibliography}
\end{document}